\documentclass[MSNbibl,number,citesort,seceqn,dvips]{arxbj}
\usepackage{upgreek,dcolumn}

\aid{0}
\volume{19}
\issue{5B}
\pubyear{2013}
\firstpage{2200}
\lastpage{2221}
\doi{10.3150/12-BEJ449} 

\makeatletter
\newtheorem{lem}{Lemma}[section]
\newtheorem{prp}{Proposition}[section]
\newremark{remark}{Remark}[section]
\newcolumntype{d}[1]{D{.}{.}{#1}}
\def\y{{\mathbf{y}}}
\def\t{{\mathbf{t}}}
\def\z{{\mathbf{z}}}
\def\j{{\mathbf{j}}}
\def\u{{\mathbf{u}}}
\def\A{{\mathbf{A}}}
\def\D{{\mathbf{D}}}
\def\C{{\mathbf{C}}}
\def\H{{\mathbf{H}}}
\def\I{{\mathbf{I}}}
\def\P{{\mathbf{P}}}
\def\Q{{\mathbf{Q}}}
\def\V{{\mathbf{V}}}
\def\W{{\mathbf{W}}}
\def\X{{\mathbf{X}}}
\def\bmu{{\boldsymbol{\mu}}}
\def\bbe{{\boldsymbol{\beta}}}
\def\bxi{{\boldsymbol{\xi}}}
\def\R{{\mathbf{R}}}
\def\s{{\mathbf{s}}}
\def\G{{\mathbf{G}}}
\def\bnabla{{\boldsymbol{\nabla}}}
\def\zero{{\boldsymbol{0}}}
\def\bep{{\boldsymbol{\epsilon}}}

\def\ep{{\varepsilon}}
\def\la{{\lambda}}
\def\lah{{\hat\la}}
\def\bbeh{{\widehat\bbe}}
\def\muh{{\widehat\mu}}
\def\bmuh{{\widehat\bmu}}
\def\Ga{{\Gamma}}
\def\bxih{{\widehat\bxi}}
\def\dd{\mathrm{d}}
\def\Nc{{\mathcal{N}}}
\def\infi{{\infty}}
\def\De{{\Delta}}
\def\Deh{{\widehat\De}}
\newcommand{\Var}{\operatorname{Var}}
\newcommand{\tr}{\operatorname{tr}}
\newcommand{\Ch}{\operatorname{Ch}}
\newcommand{\diag}{\operatorname{diag}}

\makeatother

\begin{document}
\begin{frontmatter}

\title{Dominance properties of constrained Bayes and empirical Bayes
estimators}
\runtitle{Dominance properties of constrained Bayes estimators}

\begin{aug}
\author[1]{\fnms{Tatsuya} \snm{Kubokawa}\corref{}\thanksref{1}\ead[label=e1]{tatsuya@e.u-tokyo.ac.jp}} \and
\author[2]{\fnms{William E.} \snm{Strawderman}\thanksref{2}\ead[label=e2]{straw@stat.rutgers.edu}}
\runauthor{T. Kubokawa and W.E. Strawderman} 
\address[1]{Faculty of Economics, University of Tokyo,
7-3-1, Hongo, Bunkyo-ku, Tokyo 113-0033, Japan. \\\printead{e1}}
\address[2]{Department of Statistics, Rutgers University,
Piscataway, NJ 08854-8019, USA. \\\printead{e2}}
\end{aug}

\received{\smonth{7} \syear{2011}}
\revised{\smonth{2} \syear{2012}}

%
\begin{abstract}
This paper studies decision theoretic properties of benchmarked
estimators which are of some importance in small area estimation problems.
Benchmarking is intended to improve certain aggregate properties (such
as study-wide averages) when model based estimates have been applied to
individual small areas.
We study decision-theoretic properties of such estimators by reducing
the problem to one of studying these problems in a related derived problem.
For certain such problems, we show that unconstrained solutions in the
original (unbenchmarked) problem give unconstrained Bayes and improved
estimators which automatically satisfy the benchmark constraint.
Also, dominance properties of constrained empirical Bayes estimators
are shown in the Fay--Herriot model, a frequently used model in small
area estimation.
\end{abstract}

%
\vspace*{-6pt}

\end{frontmatter}

\section{Introduction}\label{secint}

This paper studies decision theoretic properties of benchmarked
estimators which are of some importance in small area estimation problems.
Benchmarking is intended to improve certain aggregate properties (such
as study-wide averages) when empirical Bayes estimates have been
applied to individual small areas.
For example, model based small area estimates are often such that the
average of a particular estimate over all areas may differ
substantially from the average derived from a direct estimate.
The reader is referred to the articles of Datta \textit{et al.} \cite{r6} for an extended discussion of the background and
desirability of benchmarking.
Also see Frey and Cressie \cite{r10}, Ghosh \cite{r11} and Pfeffermann and Tiller \cite{r16} for related issues.
For good accounts of small area estimation, see Battese, Harter and Fuller \cite{r1}, Prasad and Rao \cite{r17}, Ghosh and Rao \cite{r12}, Rao \cite{r18}
and Datta, Rao and Smith \cite{r8}.

A useful method for benchmarking is the constrained Bayes and empirical
Bayes estimator suggested by Ghosh \cite{r11}.\vadjust{\goodbreak}
Since the constrained Bayes estimator is not a real Bayesian procedure,
its decision-theoretic properties like admissibility and minimaxity are
interesting questions, though little has been known about such properties.
Another query is whether there exists a prior distribution which
results in the (unconstrained) real Bayes estimator satisfying the constraint.
This paper will address these problems in a decision-theoretic framework.

In Section \ref{secgeneral}, we begin by explaining the empirical
Bayes estimators of small-area means and their benchmarking in the
Fay--Herriot area-level model, and give a motivation as well as the
setup of the problem.
To investigate basic decision-theoretic properties of the constrained
estimator, we decompose the risk function into two pieces; one depends
on the risk of the unconstrained estimator in a related problem and one
depends on the given means and the benchmark constraint but not the
estimator in question.
Admissibility considerations and sometimes minimaxity are then reduced
to the study of these properties in a related problem.
Section \ref{secunc} studies prior distributions in the original
problem that result in Bayes estimators which automatically satisfy the
benchmark constraint.
In fact, we clarify a condition on such prior distributions and gives examples.
Such prior distributions and the resulting Bayes estimators enable us
to study admissibility.
The results in Section \ref{secgeneral} are given without assuming
normality of the underlying distribution.

Section \ref{secnormal} assumes the multivariate normal distribution,
and provides more detailed properties for minimaxity, admissibility and
inadmissibility of constrained Bayes estimators.
In Section~\ref{secnunc}, we present a prior distribution such that
the resulting (generalized) Bayes estimator satisfies the constraint
and is also minimax.
Admissibility and minimaxity of such unconstrained Bayes estimators are
discussed based on some preliminary results given in the literature.

As indicated above, benchmarking is useful in the framework of small
area estimation.
The Fay--Herriot model is one that is often utilized in small area
estimation problems.
In Section \ref{secFH}, we consider this model and investigate
conditions under which a constrained empirical Bayes estimator improves
on the constrained uniform-prior generalized Bayes estimator, namely
the constrained direct estimator.
Since the Fay--Herriot model has heteroscedastic variances and employs
covariates as regressors, establishing minimaxity of the constrained
empirical Bayes estimator, while somewhat challenging, seems to be
potentially useful.
We also consider a prior distribution which results in an unconstrained
empirical Bayes estimator satisfying the constraint and minimaxity.
These constrained and unconstrained empirical Bayes estimators are
investigated in terms of their risk performances by simulation as well
as in terms of the conditions for their improvement or minimaxity.
Finally, some concluding remarks are given in Section \ref{secremark}.

\section{The constrained problem and the dominance property }\label{secgeneral}

\subsection{The area-level model and the setup of the problem}

The Fay--Herriot model has been used as an area-level model in
small-area estimation.
Let $y_1, \ldots, y_k$ be the direct estimators of the $k$ small-area
means $\mu_1, \ldots, \mu_k$.
The direct estimator may be taken to be a crude estimator like a sample
mean over the small area.
This is modeled as
%
\begin{equation}
\label{eqnM1} y_i = \mu_i + \ep_i,\qquad i=1,
\ldots, k,
\end{equation}
where $\ep_1, \ldots, \ep_k$ are independently distributed as $\ep_i\sim
\Nc(0, d_i)$.
While the values of $y_i$'s are reported from government agencies, the
values of the variances $d_i$'s, are usually not available, and we need
to get the values by estimation from past data or other methods.
In the framework of small area estimation, the $d_i$'s are treated as
known constants.
Small area refers to a small geographical area or a group for which
little information is obtained from the sample survey, and the direct
estimator based only on the data from a given small area is likely to
be unreliable because only a few observations are available from the
small area.
Also, $y_i$ is more unreliable for larger $d_i$.
To increase the precision of the estimate, relevant supplementary
information such as data from other related small areas or data on
covariates is used through Bayesian models.
Fay and Herriot \cite{r9} suggested a Bayesian model for $\mu_i$ in (\ref
{eqnM1}) with prior distribution of $\mu_i$ given by
%
\begin{equation}
\label{eqnM2} \mu_i \sim\Nc\bigl(\mathbf{x}_i'
\bbe, \la\bigr),\qquad i=1, \ldots, k,
\end{equation}
where $\mathbf{x}_i$ is a $p$-variate known vector including
covariates, $\bbe$
is a $p$-variate unknown vector and $\la$ is an unknown variance.
The resulting empirical Bayes estimator of $\mu_i$ is
\[
\muh_i^{\mathrm{EB}} = \mathbf{x}_i'\bbeh+
\frac{\lah}{\lah+ d_i} \bigl(y_i - \mathbf{x}_i'
\bbeh\bigr),
\]
where $\bbeh$ and $\lah$ are suitable estimators of $\bbe$ and $\la$.
For larger $d_i$, $\muh_i^{\mathrm{EB}}$ can shrink $y_i$ more toward the
estimator $\mathbf{x}_i'\bbeh$, so that it is expected that
$\muh_i^{\mathrm{EB}}$ has a
higher precision than $y_i$.
To measure the uncertainty of $\muh_i^{\mathrm{EB}}$, Prasad and Rao \cite{r17},
Datta and Lahiri \cite{r7} and Datta, Rao and Smith \cite{r8} derived a second-order
approximation of the mean squared error (MSE) of $\muh_i^{\mathrm{EB}}$ for
large $k$ under the unconditional model of (\ref{eqnM1}) and (\ref
{eqnM2}), namely, $y_i\sim\Nc(\mathbf{x}_i'\bbe, \la+d_i)$.
Since the second-order approximation of the MSE is smaller in a large
parameter space than that of the direct estimator, $\muh_i^{\mathrm{EB}}$ has
been used practically.
However, it is not guaranteed analytically that $\muh_i^{\mathrm{EB}}$ has a
uniformly smaller MSE than the direct estimator in terms of minimizing
the second-order approximation, much less the exact MSE.
This point will be demonstrated in Section \ref{secFH1}.

We can consider the uncertainty of an estimator $\muh_i$ through the
two kinds of MSE: the conditional MSE $E[(\muh_i-\mu_i)^2 | \mu_i]$
given $\mu_i$ and the unconditional MSE $E[(\muh_i-\mu_i)^2 ]$.
The unconditional MSE is measured based on the unconditional (marginal)
distribution of (\ref{eqnM1}) and (\ref{eqnM2}), and it is
interpreted as a Bayesian measure from a Bayesian perspective.
The conditional MSE is a measure supported by a frequentist, and it is
stronger since it does not assume a distribution for $\mu_i$.
In the framework of the conditional MSE, it is known that $y_i$ is
admissible in the estimation of the individual mean $\mu_i$, namely,
$\muh_i^{\mathrm{EB}}$ does not improve on $y_i$ uniformly in terms of the
conditional MSE.
However, in simultaneous estimation of the small area means $\mu_i$,
$i=1, \ldots, k$, $\muh_i^{\mathrm{EB}}$ improves on the direct estimator $y_i$
for $k\geq3$ due to the Stein effect (Stein \cite{r22}).
Thus, the framework of simultaneous estimation can justify the
improvement of $\muh_i^{\mathrm{EB}}$ theoretically.

In this paper, we consider simultaneous estimation of the small-area means.
It is convenient to handle the problem in matricial form.
Let $\y=(y_1, \ldots, y_k)'$, $\bmu=(\mu_1, \ldots, \mu_k)'$ and
$\bep
=(\ep_1, \ldots, \ep_k)'$.
Then, the model (\ref{eqnM1}) is written as
%
\begin{equation}
\label{eqnM} \y=\bmu+ \bep,
\end{equation}
where $\bep\sim\Nc_k(\zero, \D)$ for $\D=\diag(d_1, \ldots, d_k)$,
the $k\times k$ diagonal matrix.
When we estimate $\bmu$ by $\bmuh=\bmuh(\y)=(\muh_1, \ldots, \muh_k)'$
based on $\y$, the estimator is evaluated in terms of the conditional
risk function given $\bmu$,
\[
R(\bmu, \bmuh)=E \bigl[L(\bmu, \bmuh;\Q) | \bmu \bigr],
\]
relative to weighted squared error loss
%
\begin{equation}
\label{eqnloss} L(\bmu, \bmuh;\Q) =(\bmuh-\bmu)'\Q(\bmuh-\bmu),
\end{equation}
where $\Q$ is a positive definite matrix.
In a decision-theoretic framework, the set of direct estimators $\y$ is
minimax, but inadmissible by the so-called Stein effect for $k\geq3$,
namely, there exist shrinkage or empirical Bayes estimators which have
uniformly smaller risks than $\y$ for large $k$.
Since a goal in small area estimation is the derivation of estimators
having high precisions, desirable estimators should satisfy at lease
the requirement that they have uniformly smaller risks than $\y$.
This corresponds to the derivation of estimators which are minimax or
improve on $\y$ in terms of $R(\bmu, \bmuh)$.
It is noted that if an estimator improves on $\y$ in terms of the
conditional risk $R(\bmu,\bmuh)$, then it improves on $\y$ relative to
the unconditional risk
\[
R^U(\pi, \bmuh)=E^{\pi} \bigl[ E \bigl[L(\bmu, \bmuh;\Q) |
\bmu \bigr] \bigr], 
\]
where $\pi(\bmu)$ is a distribution of $\bmu$.
The unconditional risk is treated for the Fay--Herriot model in Section
\ref{secFH}.

As indicated in Section \ref{secint}, a drawback of the empirical
Bayes estimator $\muh_i^{\mathrm{EB}}$'s is that the weighted sum $\sum_{i=1}^k
w_i \muh_i^{\mathrm{EB}}$ is not equal to $\sum_{i=1}^k w_i y_i$, which, for
example, corresponds to the total sample mean over the whole area,
where $w_i$'s are nonnegative constants.
In the literature, several methods have been proposed in order to
\vspace{1pt}
benchmark an estimator $\bmuh$ so as to satisfy the constraint $\sum_{i=1}^k w_i \muh_i=\sum_{i=1}^k w_i y_i$.
Of these, Ghosh \cite{r11} suggested the constrained Bayes estimator to
satisfy the constraint.
In this paper, we consider the general constraint given by
%
\begin{equation}
\label{eqnconst0} \W'\bmuh=\t(\y),
\end{equation}
where $\W$ is a $k\times m$ matrix with rank $m$, $m<k$, and $\t=\t
(\y
)$ is a function from $\R^k$ to $\R^m$.
Typical examples of $\t(\y)$ are $t(\y)=\sum_{i=1}^k w_i y_i$ and
$t(\y
)=t_0$, a constant.
Denote the class of benchmarked estimators by
\[
\Ga_B=\bigl\{\bmuh\in\Ga\vert \W' \bmuh= \t(\y)\bigr\},
\]
where $\Ga$ is the class of estimators with second moments given by
$\Ga= \{ \bmuh| E[\bmuh'\bmuh| \bmu]<\infi\}$.
When a prior distribution $\pi$ is assumed for $\bmu$, the constrained
Bayes estimator is defined as the estimator $\bmuh$ which minimizes the
posterior risk function $E^{\pi}[ (\bmuh-\bmu)'\Q(\bmuh-\bmu) |
\y]$
subject to $\bmuh\in\Ga_B$, where $E^{\pi}[\cdot|\y]$ denotes a
posterior expectation given $\y$.
Noting that
\[
E^{\pi}\bigl[(\bmuh-\bmu)'\Q(\bmuh-\bmu) | \y
\bigr]=E^{\pi}\bigl[\bigl(\bmuh^B-\bmu\bigr)'\Q
\bigl(\bmuh^B-\bmu\bigr) | \y\bigr]+\bigl(\bmuh-\bmuh^B
\bigr)'\Q\bigl(\bmuh-\bmuh^B\bigr)
\]
for the Bayes estimator $\bmuh^B=E^{\pi}[\bmu| \y]$, Datta \textit{et al.} \cite{r6} showed that the constrained Bayes estimator is given by
\[
\bmuh^{\mathrm{CB}} = \bmuh^B +\Q^{-1}\W\bigl(
\W'\Q^{-1}\W\bigr)^{-1}\bigl\{ \t(\y)-
\W'\bmuh^B\bigr\}, 
\]
Motivated by the constrained Bayes estimator, we can construct the
following constrained estimator based on any given estimator $\bmuh$:
%
\begin{equation}
\label{eqnCE} \bmuh^C(\bmuh,\t)=\bmuh+\Q^{-1}\W\bigl(
\W'\Q^{-1}\W\bigr)^{-1}\bigl\{ \t(\y )-
\W'\bmuh \bigr\},
\end{equation}
and denote the class by
\[
\Ga_C = \bigl\{\bmuh^C(\bmuh,\t) | \bmuh\in\Ga\bigr\}.
\]
It is seen that
\[
\Ga_C \subset\Ga_B \subset\Ga.
\]
Since $\y$ is the generalized Bayes estimator of $\bmu$ against the
uniform prior, the constrained generalized Bayes estimator against the
uniform prior is expressed as
%
\begin{equation}
\label{eqnCM} \bmuh^{\mathrm{CM}}(\t)=\y+ \Q^{-1}\W\bigl(
\W'\Q^{-1}\W\bigr)^{-1}\bigl\{ \t(\y)-
\W'\y\bigr\}.
\end{equation}
It is noted that the direct estimator $\y$ satisfies the constraint
when the constraint is that $\t(\y)=\W'\y$.

Since the constrained Bayes estimator is not necessarily the Bayes
estimator among all estimators in $\Ga$, we have several interesting
questions from a decision-theoretic perspective.
For example, are the properties of minimaxity and inadmissibility of
$\y
$ inherited by the constrained estimator $\bmuh^{\mathrm{CM}}$?
Can one construct an empirical Bayes estimator improving on $\y$ or
$\bmuh^{\mathrm{CM}}$?
Such issues have not been studied in the literature to our knowledge.
The aim of this paper is to investigate such decision-theoretic
properties for the constrained estimators.

\subsection{Basic properties of a constrained estimator}\label{secCB}

In this subsection, we investigate basic properties of minimaxity and
admissibility of the constrained estimator under the constraint (\ref
{eqnconst0}) in the model (\ref{eqnM}), where normality of $\bep$ is
not assumed in this and the next subsections.
We begin by decomposing the risk function, which will be useful for
investigating the basic properties.
Let
\[
\P_{\W}=\Q^{-1}\W\bigl(\W'\Q^{-1}\W
\bigr)^{-1}\W'.
\]
Then the constrained estimator (\ref{eqnCE}) is expressed as
%
\begin{equation}
\label{eqnCE1} \bmuh^C(\bmuh,\t)=(\I-\P_{\W}) \bmuh+
\Q^{-1}\W\bigl(\W'\Q^{-1}\W \bigr)^{-1}
\t(\y).
\end{equation}
To evaluate the risk of $\bmuh^C(\bmuh,\t)$, note that $\W'(\I-\P_{\W}
)=\zero$,
\begin{eqnarray*}
\bmuh^C(\bmuh,\t)-\bmu&=&(\I-\P_{\W}) (\bmuh-\bmu) +
\Q^{-1}\W\bigl(\W'\Q^{-1}\W \bigr)^{-1}
\bigl\{ \t(\y)-\W'\bmu\bigr\},
\\
(\I-\P_{\W})'\Q(\I-\P_{\W})&=&\Q- \W\bigl(
\W'\Q^{-1}\W\bigr)^{-1}\W'=\Q (\I-
\P_{\W}).
\end{eqnarray*}
Then the conditional risk function of $\bmuh^C(\bmuh,\t)$ relative to
the loss (\ref{eqnloss}) can be decomposed into two parts as given in
the following lemma.

\begin{lem}\label{lem1}
Assume that $\bmuh\in\Ga$.
It follows that the conditional risk function of $\bmuh^C(\bmuh,\t)$
relative to the loss $L(\bmu, \bmuh;\Q)$ is expressed as
%
\begin{equation}
\label{eqnrisk0} R\bigl(\bmu, \bmuh^C(\bmuh,\t)\bigr)=R_1(
\bmu, \bmuh) + R_2(\bmu, \t),
\end{equation}
where $R_1(\bmu, \bmuh)= E [(\bmuh-\bmu)'\Q(\I-\P_{\W}
)(\bmuh-\bmu
) | \bmu ]$ and
\[
R_2(\bmu,\t) = E \bigl[\bigl(\t(\y)-\W'\bmu
\bigr)'\bigl(\W'\Q^{-1}\W\bigr)^{-1}
\bigl(\t (\y)-\W'\bmu\bigr) | \bmu \bigr].
\]
\end{lem}

Since $\t(\y)$ is a given function and $R_2(\bmu,\t)$ does not depend
on the estimator $\bmuh$, the problem of finding improved estimators
(in the original benchmark problem) can be reduced to that of finding
superior estimators $\bmuh$ in terms of the risk function $R_1(\bmu
,\bmuh)$ relative to the loss function $L(\bmu, \bmuh;\Q(\I-\P_{\W}))$.

\begin{prp}\label{prpdom}
For two estimators $\bmuh_1$ and $\bmuh_2$ in $\Ga$, and the
corresponding constrained estimators $\bmuh^C(\bmuh_1,\t)$ and
$\bmuh^C(\bmuh_2,\t)$ in $\Ga_C$,
$\bmuh^C(\bmuh_1,\t)$ dominates $\bmuh^C(\bmuh_2,\t)$ relative to the
loss $L(\bmu, \bmuh;\Q)$ if and only if $\bmuh_1$ dominates $\bmuh_2$
relative to the loss $L(\bmu, \bmuh;\Q(\I-\P_{\W}))$.
\end{prp}

This proposition implies the following proposition concerning admissibility.

\begin{prp}\label{prpadm}
Assume that $\bmuh\in\Ga$.
Then the constrained estimator $\bmuh^C(\bmuh,\t)$ is admissible in
$\Ga_C$ in terms of the risk $R(\bmu,\bmuh^C)$ if and only if $\bmuh$ is
admissible in $\Ga$ in terms of the risk $R_1(\bmu, \bmuh)$.
\end{prp}

The above propositions show that dominance properties and admissibility
of a constrained estimator $\bmuh^C(\bmuh,\t)$ can be reduced to those
of the estimator $\bmuh$ in terms of the risk $R_1(\bmu,\bmuh)$.

Concerning minimaxity, on the other hand, it is seen that the estimator
$\bmuh^C(\bmuh^*,\t)$ is minimax within the class $\Ga_C$ if and only
if $\inf_{\bmuh\in\Ga_C} \sup_{\bmu} R(\bmu, \bmuh)=\sup_{\bmu}
R(\bmu,
\bmuh^*)$, or
\[
\inf_{\bmuh\in\Ga} \sup_{\bmu} \bigl\{ R_1(\bmu, \bmuh) +
R_2(\bmu, \t) \bigr\} = \sup_{\bmu} \bigl\{ R_1
\bigl(\bmu, \bmuh^*\bigr) + R_2(\bmu, \t) \bigr\}.
\]
This condition is satisfied if there exists a sequence of prior
distributions $\{ \pi_n(\bmu) \}_{n=1,2,\ldots}$ such that
%
\begin{equation}
\label{eqnminc} \sup_{\bmu} \bigl\{ R_1\bigl(\bmu, \bmuh^*
\bigr) + R_2(\bmu, \t) \bigr\} \leq\operatorname{\lim\sup}\limits
_{n\to\infi} \int \bigl\{ R_1(\bmu, \bmuh_n)
+ R_2(\bmu, \t) \bigr\}\pi_n(\bmu)\,\dd\bmu,
\end{equation}
where $\bmuh_n$ is the corresponding Bayes estimator relative to the
risk $\int R(\bmu, \bmuh)\pi_n(\bmu)\,\dd\bmu$.
This condition follows from Theorem 6.5.2 in Zacks \cite{r23}.

\begin{prp}\label{prpmin}
If $\bmuh^*$ satisfies the condition $(\ref{eqnminc})$, then the
estimator $\bmuh^C(\bmuh^*,\t)$ is minimax within $\Ga_C$.
\end{prp}

In particular, under the following condition, the minimaxity problem
for the conditional risk $R(\bmu, \bmuh^C)$ reduces to that of the risk
$R_1(\bmu, \bmuh)$.
\begin{enumerate}[(A1)]
\item[(A1)] Assume that $R_2(\bmu, \bmuh)$ does not depend on the
unknown $\bmu$.
\end{enumerate}
%
\begin{prp}\label{prpmins}
Assume the condition {\rm(A1)}.
Then the constrained estimator $\bmuh^C(\bmuh^*,\t)$ is minimax within
$\Ga_C$ if and only if $\bmuh^*$ is minimax in terms of the risk
$R_1(\bmu,\bmuh)$ in $\Ga$.
\end{prp}

Condition (A1) is satisfied for two typical examples of $\t(\y)$:

\textit{Case 1}: $\t(\y)=\W'\y$.
In this case, it typically happens that $R_2(\bmu, \t)$ is independent
of $\bmu$ under the distributional assumption of a location family, and
the condition (A1) holds.

\textit{Case 2}: $\t(\y)=\t_0$, a constant.
In this case, we need to restrict the space of $\bmu$ to $\{ \bmu| \W'\bmu=\t_0\}$.
Then it is clear that $R_2(\bmu,\t_0)=0$ on the restricted space.

\subsection{Unconstrained Bayes estimators satisfying the constraint}\label{secunc}

In the previous subsections, we studied shrinkage estimators induced
from the constrained Bayes estimator and investigated their
decision-theoretic properties within the class of constrained estimators.
In some cases, however, we can derive constrained Bayes estimators
without direct consideration of the constraint.
In this subsection, we find a condition on prior distributions such
that the resulting unconstrained generalized Bayes estimators satisfy
the constraint automatically, where the normality of $\bep$ is not assumed.

Assume a prior distribution $\pi$ for $\bmu$.
According to the expression in (\ref{eqnCE1}), we decompose $\bmu$ as
$\bmu=(\I-\P_{\W})\bmu+\P_{\W}\bmu$, which implies that the Bayes estimator
of $\bmuh^B$ can be expressed as
\[
\bmuh^B = (\I-\P_{\W})E^{\pi}[\bmu| \y] +
\P_{\W} E^{\pi}[\bmu|\y].
\]
Comparing this expression and the constrained estimator (\ref
{eqnCE1}), we can see that the unconstrained Bayes estimator $\bmuh^B$
belongs to the class $\Ga_C$ if the prior distribution satisfies the equation
\[
\P_{\W} E^{\pi}[\bmu|\y] = \Q^{-1}\W\bigl(
\W'\Q^{-1}\W\bigr)^{-1}\t(\y). 
\]
It follows from the definition of $\P_{\W}$ that this equality is
simplified as
%
\begin{equation}
\label{eqnunc2} \W' E^{\pi}[\bmu|\y] = \t(\y).
\end{equation}

Since the condition (\ref{eqnunc2}) means that the posterior
expectation of $\W'\bmu$ is $\t(\y)$, the following transformation is
convenient for investigating prior distributions satisfying (\ref{eqnunc2}).
Let $\H$ be a $k\times k$ orthogonal matrix such that
%
\begin{equation}
\label{eqnH} \H\Q^{-1/2}\W\bigl(\W'\Q^{-1}\W
\bigr)^{-1}\W'\Q^{-1/2}\H'=\pmatrix{
\zero_{k-m} & \zero
\cr
\zero&\I_m }.
\end{equation}
Let $\H'=(\H_1',\H_2')$ for the $k\times(k-m)$ matrix $\H_1$.
Also, let $\bxi=\H\Q^{1/2}\bmu$ and $\bxi_i=\H_i\Q^{1/2}\bmu$
for $i=1, 2$.
Then, $\P_{\W}=\Q^{-1/2}\H_2'\H_2\Q^{1/2}$ and $\I-\P_{\W}=\Q^{-1/2}\H_1'\H_1\Q^{1/2}$.
It is noted that $\H_1 \Q^{-1/2}\W=\zero$, since $\H_1 \Q^{-1/2}\W=\H_1\H_1'\H_1 \Q^{-1/2}\W=\Q^{1/2}(\I-\P_{\W})\Q^{-1}\W=\zero$.
Thus, $\bmu$ and the constrained estimator (\ref{eqnCE1}) are written as
%
\begin{eqnarray}
\label{eqnunc3} \bmu&=& \Q^{-1/2} \H\H'\Q^{1/2}
\bmu=\Q^{-1/2}\H_1'\bxi_1 +
\Q^{-1/2}\H_2'\bxi_2,
\nonumber
\\[-8pt]
\\[-8pt]
\bmuh^{C}(\bmuh, \t) &=& \Q^{-1/2}\H_1'
\H_1\Q^{1/2} \bmuh+ \Q^{-1}\W\bigl(
\W'\Q^{-1}\W\bigr)^{-1}\t(\y),
\nonumber
\end{eqnarray}
which shows that the unconstrained Bayes estimator belongs to $\Ga_C$ if
%
\begin{equation}
\label{eqnunc4} \H_2'E^{\pi}[
\bxi_2|\y] = \Q^{-1/2}\W\bigl(\W'
\Q^{-1}\W\bigr)^{-1}\t(\y).
\end{equation}
Noting that $\Q^{-1/2}\W=(\H_1'\H_1+\H_2'\H_2)\Q^{-1/2}\W=\H_2'\H_2\Q^{-1/2}\W$, we see that the equation~(\ref{eqnunc4}) holds if
\[
E^{\pi}[\bxi_2 | \y] = \H_2\Q^{-1/2}\W
\bigl(\W'\Q^{-1}\W\bigr)^{-1}\t (\y).
\]
For example, consider the case that the constraint is given by $\t(\y
)=\W'\Q^{-1/2}\s(\y)$ for a $k$-variate vector $\s(\y)$ of
functions of
$\y$.
In this case, we have that
\[
\H_2\Q^{-1/2}\W\bigl(\W'\Q^{-1}\W
\bigr)^{-1}\W'\Q^{-1/2}\H'\H\s(\y )=
\H_2\s(\y),
\]
so that the condition (\ref{eqnunc4}) may be simplified as $E^{\pi
}[\bxi_2 | \y] = \H_2 \s(\y)$.
Thus, we summarize the condition in the following.
\begin{enumerate}[(A2)]
\item[(A2)] Assume that $\W' E^{\pi}[\bmu| \y] =\t(\y)$, or that
$E^{\pi}[\bxi_2 | \y] = \H_2\Q^{-1/2}\W(\W'\Q^{-1}\W)^{-1}\t
(\y)$.
The latter condition is simplified as $E^{\pi}[\bxi_2 | \y] = \H_2
\s(\y
)$ when $\t(\y)=\W'\Q^{-1/2}\s(\y)$.
\end{enumerate}
\begin{prp}\label{prpunc}
The unconstrained Bayes estimators belong to the class $\Ga_C$, namely
they automatically satisfy the constraint $\W'\bmuh=\t(\y)$ if the
posterior expectation $E^{\pi}[\bmu|\y]$ or $E^{\pi}[\bxi_2 | \y]$
satisfies the condition {\rm(A2)}.
\end{prp}

\textit{Case 1}: $\t(\y)=\W'\y$.
In this case, the condition (A2) is $\W' E^{\pi}[\bmu| \y] =\W'\y
$ or
$E^{\pi}[\bxi_2 | \y] = \H_2 \Q^{1/2} \y$ for $\bxi_2=\H_2\Q^{1/2}\bmu$.
As explained in the next section, it suffices that we assume the
uniform prior for $\bxi_2$ under normality of $\bep$.

\textit{Case 2}: $\t(\y)=\t_0$, a constant.
In this case, the condition (A2) is $E^{\pi}[\bxi_2 | \y] = \H_2\Q^{-1/2}\W\*(\W'\Q^{-1}\W)^{-1}\t_0$, which suggests that $\bxi_2$ should
take a point mass at $\bxi_2=\H_2\Q^{-1/2}\W\*(\W'\Q^{-1}\W
)^{-1}\t_0$.
Since $\W'\Q^{-1/2}\H_1'=\zero$, it is verified that this restriction
satisfies $\W'\bmu=\t_0$.

\section{Properties under normality and conditional risk}\label{secnormal}

In this section, we further investigate minimaxity and admissibility
properties for the benchmark problem in the model (\ref{eqnM}), where
normality of $\bep$ is assumed.

\subsection{Constrained Bayes estimator}

We begin by deriving the canonical form of the model (\ref{eqnM}) with
$\y$ having a multivariate normal distribution $\Nc_k(\bmu, \D)$.
For the matrix $\H$ defined by (\ref{eqnH}), let $\z_i=\H_i\Q^{1/2}\y$
and $\V_{ij}=\H_i\Q^{1/2}\D\Q^{1/2}\H_j'$ for $i, j= 1, 2$.
Then, $\z=(\z_1', \z_2')'$ is distributed as
%
\begin{equation}
\label{eqncan2} %
\pmatrix{\z_1
\cr
\z_2 } \sim
\Nc_k \left( %
\pmatrix{ \bxi_1
\cr
\bxi_2 %
}, \pmatrix{ %
\V_{11}&
\V_{12}
\cr
\V_{21}&\V_{22} %
}\right ),
\end{equation}
where $\bxi_i=\H_i\Q^{1/2}\bmu$, $i=1, 2$, for $\H_i$ defined in
(\ref{eqnH}).
The problem of finding a constrained Bayes estimator may be expressed
as the minimization of $E^{\pi}[(\bxih_1-\bxi_1)'(\bxih_1-\bxi_1)|\z]$
subject to $\W'\Q^{-1/2}\H_2'\bxih_2=\t(\y)$, since $\W'\Q^{-1/2}\H_1'=\zero$.
The constrained estimators given in (\ref{eqnCE}) and (\ref{eqnCM})
are rewritten as
\begin{eqnarray*}
\label{eqncan3} %
\bmuh^C(\bmuh,\t)&=& \Q^{-1/2}
\H_1'\bxih_1 + \Q^{-1}\W\bigl(
\W'\Q^{-1}\W\bigr)^{-1}\t(\y) \equiv
\bmuh^{C*}(\bxih_1,\t),
\\
\bmuh^{\mathrm{CM}}(\t) &=& \Q^{-1/2}\H_1'
\z_1 + \Q^{-1}\W\bigl(\W'\Q^{-1}\W
\bigr)^{-1}\t(\y) \equiv\bmuh^{C*}(\z_1,\t).
\end{eqnarray*}
For $\bxih_1=\H_1\Q^{1/2}\bmuh$, and the conditional risk $R_1(\bmu
,\bmuh)$ given in (\ref{eqnrisk0}) is written as
\[
R_1(\bmu,\bmuh)= E\bigl[\Vert\bxih_1-
\bxi_1 \Vert^2|\bxi_1\bigr]=R^*(
\bxi_1,\bxih_1), 
\]
where $\Vert\bxih_1-\bxi_1 \Vert^2=(\bxih_1-\bxi_1)'(\bxih_1-\bxi_1)$.
Hence from Proposition \ref{prpadm}, we get the following proposition.

\begin{prp}\label{prpadnor}
If $\bxih_1$ is admissible in terms of $R^*(\bxi_1,\bxih_1)$, then
$\bmuh$ is admissible\vspace*{2pt} within the class $\Ga_C$.
In particular, if $\bxih_1$ is the Bayes estimator for a proper prior
on $\bxi_1$, then $\bmuh$ is admissible within $\Ga_C$.
If $\bxih_1$ is inadmissible in terms of the risk $R^*(\bxi_1,\bxih_1)$, then $\bmuh$ is inadmissible.
\end{prp}

Also, from Propositions \ref{prpmins} and \ref{prpadnor} and the
well-known results of James and Stein \cite{r13} and Brown \cite{r5}, the next
proposition follows.
%
\begin{prp}\label{prpdec}
The constrained generalized Bayes estimator $\bmuh^{\mathrm{CM}}(\t)$ for the
uniform prior has the following decision-theoretic properties:
\begin{enumerate}[(2)]

\item [(1)] $\bmuh^{\mathrm{CM}}(\t)$ is minimax within $\Ga_C$ under the condition
{\rm(A1)}.

\item [(2)] $\bmuh^{\mathrm{CM}}(\t)$ is admissible within $\Ga_C$ when $k-m$ is one
or two.

\item [(3)] $\bmuh^{\mathrm{CM}}(\t)$ is inadmissible within $\Ga_C$ when $k-m \geq3$.
\end{enumerate}
\end{prp}

Proposition \ref{prpdec}(3) implies that there exist shrinkage
estimators like empirical Bayes estimators which improve on $\bmuh^{\mathrm{CM}}(\t)$ for large $k$.
One of such improved estimators is given in Section \ref{secFH}.
Noting that $\z_1\sim\Nc_{k-m}(\bxi_1, \V_{11})$, from the result in
Berger \cite{r2}, we can get an admissible and minimax estimator, denoted
by $\bxih_1^{\mathrm{GB}}(\z_1,\V_{11})$, based on $(\z_1, \V_{11})$
relative to
the risk $R^*(\bxi_1, \bxih_1)$.
This leads to the constrained generalized Bayes estimator
\[
\bmuh^C\bigl(\bmuh^{\mathrm{GB}},\t\bigr)=\Q^{-1/2}
\H_1'\bxih_1^{\mathrm{GB}}(
\z_1, \V_{11}) + \Q^{-1}\W\bigl(\W'
\Q^{-1}\W\bigr)^{-1}\t(\y),
\]
which is admissible within the class $\Ga_C$ and improves on $\bmuh^{\mathrm{CM}}(\t)$ when $k-m\geq3$.

\subsection{Unconstrained Bayes estimators}\label{secnunc}

We now construct unconstrained Bayes estimators satisfying the
constraint automatically in the two cases $\t(\y)=\W'\y$ and $\t
(\y)=\t_0$, a constant.
To this end, the following decomposition is useful:
%
\begin{equation}
\label{eqndec2} \pmatrix{ %
\z_3
\cr
\z_2
}\sim\Nc \left( \pmatrix{ %
\bxi_3
\cr
\bxi_2 %
}, %
\pmatrix{\V_{11.2}&\zero
\cr
\zero& \V_{22} } \right),
\end{equation}
where $\z_3=\z_1-\V_{12}\V_{22}^{-1}\z_2$, $\bxi_3=\bxi_1-\V_{12}\V_{22}^{-1}\bxi_2$ and $\V_{11.2}=\V_{11}-\V_{12}\V_{22}^{-1}\V_{21}$.
Note that $\z_1=\z_3 +\V_{12}\V_{22}^{-1}\z_2$ and that $\z_3$ is
independent of $\V_{12}\V_{22}^{-1}\z_2$.

\textit{Case 1}: $\t(\y)=\W'\y$.
Consider the decomposition (\ref{eqndec2}).
Note that from (\ref{eqnunc3}), $\bmu$ is written as
\[
\bmu= \Q^{-1/2}\H_1' \bxi_3 +
\Q^{-1/2}\bigl(\H_2'+\H_1'
\V_{12}\V_{22}^{-1}\bigr)\bxi_2.
\]
Assume a prior distribution $\pi(\bxi_3)$ for $\bxi_3$ and the uniform
prior $\pi(\bxi_2)=1$ for $\bxi_2$.
Then, the resulting generalized Bayes estimator is
%
\begin{eqnarray}
\label{eqnGB1} \bmuh^{\mathrm{GB}1} &=& \Q^{-1/2}\H_1'
\bxih^{\mathrm{GB}}_3(\z_3,\V_{11.2}) +
\Q^{-1/2}\bigl(\H_2'+\H_1'
\V_{12}\V_{22}^{-1}\bigr)\z_2
\nonumber
\\[-8pt]
\\[-8pt]
&=& \Q^{-1/2}\H_1' \bigl\{
\bxih^{\mathrm{GB}}_3(\z_3,\V_{11.2}) +
\V_{12}\V_{22}^{-1}\z_2\bigr\} +
\Q^{-1/2}\H_2'\z_2,
\nonumber
\end{eqnarray}
where $\bxih^{\mathrm{GB}}_3=\bxih^{\mathrm{GB}}_3(\z_3,\V_{11.2})$ is the generalized
Bayes estimator of $\bxi_3$ which can be constructed via the model $\z_3|\bxi_3\sim\Nc_{k-m}(\bxi_3,\V_{11.2})$ and $\bxi_3\sim\pi
(\bxi_3)$.
Note that the conditional risk $R(\bmu,\bmuh^{\mathrm{GB}1})=E[(\bmuh^{\mathrm{GB}1} -
\bmu)'\Q(\bmuh^{\mathrm{GB}1} - \bmu)|\bmu]$ is evaluated as
%
\begin{eqnarray}
\label{eqnrGB1} &&E\bigl[\bigl\{\bxih^{\mathrm{GB}}_3-
\bxi_3 + \V_{12}\V_{22}^{-1}(
\z_2-\bxi_2)\bigr\}'\H_1
\H_1'\bigl\{\bxih^{\mathrm{GB}}_3-
\bxi_3 + \V_{12}\V_{22}^{-1}(
\z_2-\bxi_2)\bigr\}\bigr]
\nonumber
\\
&&\qquad{}+ E\bigl[(\z_2-\bxi_2)'
\H_2\H_2'(\z_2-
\bxi_2)\bigr]
\nonumber
\\
&&\quad= E\bigl[\bigl\Vert\bxih^{\mathrm{GB}}_3-\bxi_3
\bigr\Vert^2\bigr] + E\bigl[(\z_2-\bxi_2)'
\V_{22}^{-1}\V_{21}\V_{12}
\V_{22}^{-1}(\z_2-\bxi_2)\bigr]
\\
&&\qquad{}+ E\bigl[\Vert\z_2-\bxi_2\Vert^2
\bigr]
\nonumber
\\
&&\quad= R^*\bigl(\bxi_3, \bxih^{\mathrm{GB}}_3\bigr)
+ \tr\bigl[\V_{12}\V_{22}^{-1}\V_{21}
\bigr] + \tr [\V_{22}],
\nonumber
\end{eqnarray}
where $R^*(\bxi_3, \bxih_3)=E[\Vert\bxih_3-\bxi_3\Vert^2|\bxi_3]$.
If $R^*(\bxi_3, \bxih^{\mathrm{GB}}_3)\leq R^*(\bxi_3, \z_3)$, then $R(\bmu
,\bmuh^{\mathrm{GB}1})\leq\tr[\V_{11}+\V_{22}]=\tr[\D\Q]$, since $R^*(\bxi_3,
\z_3)=\tr[\V_{11.2}]$.
Since $\tr[\D\Q]$ is the minimax risk, the unconstrained Bayes
estimator $\bmuh^{\mathrm{GB}1}$ is minimax in $\Ga$.
Noting that $\z_3\sim\Nc_{k-m}(\bxi_3, \V_{11.2})$, from the
result in
Berger \cite{r2}, we can get an admissible and minimax estimator based on
$(\z_3, \V_{11.2})$ relative to the risk $R^*(\bxi_3, \bxih_3)$.

\begin{prp}\label{prpprior1}
Assume the uniform prior for $\bxi_2$.
Then the generalized Bayes estimator $\bmuh^{\mathrm{GB}1}$ satisfies the
constraint, namely $\bmuh^{\mathrm{GB}1}\in\Ga_C$.
\begin{enumerate}[(2)]
\item [(1)] If $R^*(\bxi_3, \bxih^{\mathrm{GB}}_3)\leq R^*(\bxi_3, \z_3)$, then
$\bmuh^{\mathrm{GB}1}$ is minimax in $\Ga$.

\item [(2)] If $\bxih^{\mathrm{GB}}_3$ is admissible in terms of the risk
$R^*(\bxi_3, \cdot)$, then $\bmuh^{\mathrm{GB}1}$ is admissible within the
constrained class $\Ga_C$.

\item [(3)] When $m\geq3$, $\bmuh^{\mathrm{GB}1}$ is not admissible in the
unconstrained problem even if $\bxih^{\mathrm{GB}}_3$ is admissible in terms of
the risk $R^*(\bxi_3, \cdot)$.
\end{enumerate}
\end{prp}

\textit{Case 2}: $\t(\y)=\t_0$, a constant.
Assume that $\W'\Q^{-1/2}\H_2'$ is non-singular.
Since $\W'\bmu=\W'\Q^{-1/2}\H_2'\bxi_2=\t_0$, we can define
$\bxi_0$ by
$\bxi_0=(\W'\Q^{-1/2}\H_2')^{-1}\t_0$.
Let $\z_4=\z_1-\V_{12}\V_{22}^{-1}(\z_2-\bxi_0)$ and $\bxi_4=\bxi_1-\V_{12}\V_{22}^{-1}(\bxi_2-\bxi_0)$.
Then from the decomposition (\ref{eqndec2}), the joint distribution of
$(\z_4, \z_2)$ follows that
%
\begin{equation}
\label{eqndec3} \pmatrix{ %
\z_4
\cr
\z_2
}\sim\Nc \left(\pmatrix{ %
\bxi_4
\cr
\bxi_2 %
},\pmatrix{ %
\V_{11.2}&\zero
\cr
\zero& \V_{22} %
} \right).
\end{equation}
Assume a prior distribution $\pi(\bxi_4)$ for $\bxi_4$ and $P^{\pi
}[\bxi_2=\bxi_0]=1$ for $\bxi_2$.
Then from (\ref{eqnGB1}), the resulting generalized Bayes estimator is
%
\begin{equation}
\label{eqnGB2} \bmuh^{\mathrm{GB}2} = \Q^{-1/2}\H_1'
\bigl\{\bxih^{\mathrm{GB}}_4(\z_4,\V_{11.2}) +
\V_{12}\V_{22}^{-1}\bxi_0\bigr\} +
\Q^{-1/2}\H_2'\bxi_0,
\end{equation}
where $\bxih^{\mathrm{GB}}_4=\bxih^{\mathrm{GB}}_4(\z_4,\V_{11.2})$ is the generalized
Bayes estimator of $\bxi_4$ which can be constructed via the model $\z_4|\bxi_4\sim\Nc_{k-m}(\bxi_4,\V_{11.2})$ and $\bxi_4\sim\pi
(\bxi_4)$.
It also follows from (\ref{eqnrGB1}) that the risk function of $\bmuh^{\mathrm{GB}2}$ is
\[
R\bigl(\bmu,\bmuh^{\mathrm{GB}2}\bigr) = R^*\bigl(\bxi_4,
\bxih^{\mathrm{GB}}_4\bigr) + (\bxi_0-
\bxi_2)'\bigl\{ \V_{22}^{-1}
\V_{21}\V_{12}\V_{22}^{-1} +
\I_m\bigr\} (\bxi_0-\bxi_2),
\]
so that the admissibility of $\bmuh^{\mathrm{GB}2}$ is inherited from that of
$\bxih_4^{\mathrm{GB}}$.
If the space of $\bmu$ is restricted to $\{\bmu| \W'\bmu=\t_0\}$, and
if $R^*(\bxi_4, \bxih^{\mathrm{GB}}_4)\leq R^*(\bxi_4, \z_4)$, then $R(\bmu
,\bmuh^{\mathrm{GB}2})\leq\tr[\V_{11.2}]$.
Since $\tr[\V_{11.2}]$ is the minimax risk under the restriction, the
unconstrained Bayes estimator $\bmuh^{\mathrm{GB}2}$ is minimax in $\Ga$ when
$\bmu$ is restricted.

\begin{prp}\label{prpprior2}
Assume the point mass prior for $\bxi_2$.
Then the generalized Bayes estimator\vspace*{1pt} $\bmuh^{\mathrm{GB}2}$ satisfies the
constraint, namely $\bmuh^{\mathrm{GB}2}\in\Ga_C$.
\begin{enumerate}[(2)]
\item [(1)] If the estimator $\bxih_4^{\mathrm{GB}}$ of $\bxi_4$ is admissible
in terms of the risk $R^*(\bxi_4, \bxih^{\mathrm{GB}}_4)$, then $\bmuh^{\mathrm{GB}2}$ is
admissible in $\Ga$ $($and also $\Ga_C)$.

\item [(2)] If $R^*(\bxi_4, \bxih^{\mathrm{GB}}_4)\leq R^*(\bxi_4, \z_4)$, then
$\bmuh^{\mathrm{GB}2}$ is minimax within the class $\Ga_C$.
Further, it is minimax in $\Ga$ when $\bmu$ is restricted to $\W'\bmu=\t_0$ or $\bxi_2=\bxi_0$.
\end{enumerate}
\end{prp}

\section{Benchmarking in the Fay--Herriot model}\label{secFH}

As mentioned in the introduction and as explained in Datta \textit{et al.} \cite{r6}
benchmarking is useful in the framework of small area estimation. The
Fay--Herriot model is often utilized in such problems. In this section,
we develop a constrained empirical Bayes estimator for this model and
investigate the dominance properties.

\subsection{Constrained empirical Bayes estimator}\label{secFH1}

The Fay--Herriot model given in (\ref{eqnM1}) and (\ref{eqnM2}) can be
described in matricial form as
\[
\y| \bmu\sim \Nc_k(\bmu, \D),\qquad \D=
\diag(d_1, \ldots, d_k),\qquad \bmu\sim\Nc_k(
\X\bbe, \la\I),
\]
where $\X=(\mathbf{x}_1, \ldots, {\mathbf{x}}_k)'$ is a
$k\times p$ matrix of explanatory
variables with rank $p$, $\bbe$ is a $p\times1$ unknown vector of
regression coefficients and $\la$ is an unknown scalar.
Suppose that $d_1\geq\cdots\geq d_k$ without any loss of generality.
Consider estimation of $\bmu$ in terms of the conditional risk $R(\bmu,
\bmuh)=E[(\bmuh- \bmu)'\Q(\bmuh- \bmu) | \bmu]$ and the unconditional
risk $R^U(\pi, \bmuh)=E[(\bmuh- \bmu)'\Q(\bmuh- \bmu)]$ where
$\pi$
denotes the distribution of $\bmu$.
The Bayes estimator (under the assumption of known $ \bbe$ and $ \la
$) is given by
\[
\bmuh^B= \X\bbe+ (\D/\la+\I)^{-1}(\y- \X\bbe)=\y- \D(\D+
\la\I )^{-1}(\y-\X\bbe).
\]
For estimation of $\la$, several estimators are known including the
Prasad--Rao estimator given by Prasad and Rao \cite{r17}, the Fay--Herriot
estimator suggested by Fay and Herriot \cite{r9}, the maximum likelihood
estimator (MLE) and the restricted maximum likelihood estimator (REML).
For the MLE and REML, see Searle, Casella and McCulloch \cite{r19} and
Kubokawa \cite{r14}, for example.
Denoting an estimator of $\la$ by $\lah$, we get the empirical Bayes
estimator $
\bmuh^{\mathrm{EB}}(\lah)= \y- \D(\D+ \lah\I)^{-1}(\y-\X\bbeh(\lah
))$, where
$\bbeh(\lah)= \{\X'\V(\lah)^{-1}\X\}^{-1}\X'\V(\lah)^{-1}\y$
for $\V(\la
)=\D+\la\I$.
The empirical Bayes estimator is called the empirical best linear
unbiased predictor (EBLUP) in the framework of the linear mixed model,
namely the unconditional model (\ref{eqnEB}).
Define $\A(\la)$ by
\[
\A(\la) = \V(\la)^{-1} - \V(\la)^{-1}\X\bigl(
\X'\V(\la)^{-1}\X \bigr)^{-1}\X'\V(
\la)^{-1}.
\]
Then, the empirical Bayes estimator can be rewritten as
%
\begin{equation}
\label{eqnEB} \bmuh^{\mathrm{EB}}(\lah)=\y- \D\A(\lah)\y.
\end{equation}
Now consider the benchmark constraint $\W'\bmuh=\t(\y)$.
The constrained empirical Bayes estimator (CEB) based on $\bmuh^{\mathrm{EB}}(\lah)$ (as constructed in \ref{eqnCE1}) is given by
%
\begin{equation}
\label{eqnCEB} \bmuh^{\mathrm{CEB}}(\lah,\t) = (\I-\P_{\W})
\bmuh^{\mathrm{EB}}(\lah) + \Q^{-1}\W \bigl(\W'
\Q^{-1}\W\bigr)^{-1}\t(\y).
\end{equation}
Concerning the estimation of $\la$, we here treat the Fay--Herriot
estimator $\lah$ given by $\lah=\max\{\la^*, 0\}$ where $\la^*$ is the
solution of the equation
%
\begin{equation}
\label{eqnla} \y'\A\bigl(\la^*\bigr) \y= k-p.
\end{equation}

\textit{1. Conditional risk.} A sufficient condition for $\bmuh^{\mathrm{CEB}}(\lah,\t)$ to improve on $\bmuh^{\mathrm{CM}}(\t)$ in terms of the
conditional risk is given in the following proposition which will be
proved in the \hyperref[app]{Appendix}.

\begin{prp}\label{prpFHmin}
The constrained empirical Bayes estimator $\bmuh^{\mathrm{CEB}}(\lah,\t)$ with
$\lah$ given in $(\ref{eqnla})$ improves on $\bmuh^{\mathrm{CM}}(\t)$ given in
$(\ref{eqnCM})$ in terms of the conditional risk if the following
inequality holds:
%
\begin{equation}
\label{eqnFHcond} \min_{\la>0} \biggl\{ \frac{\tr[\D\Q_{\W}\D\A(\la)]}{\Ch_{\mathrm{max}}(\D\Q_{\W}\D
\A(\la)) } \biggr\} \geq
\frac{k-p}{2} + 2,
\end{equation}
where $\Q_{\W}=\Q-\W(\W'\Q^{-1}\W)^{-1}\W'$, and $\Ch_{\mathrm{max}}(\C
)$ denotes
the maximum eigenvalue of the matrix $\C$.
If the constraint is given by $\t(\y)=\W'\y$, then the estimator
$\bmuh^{\mathrm{CEB}}(\lah,\t)$ is minimax under the condition $(\ref{eqnFHcond})$.
\end{prp}

To derive explicit sufficient conditions, it is noted that
\begin{eqnarray*}
\tr\bigl[\D\Q_{\W}\D\A(\la)\bigr]&=& \tr\bigl[\D\Q_{\W}\D(
\D+\la\I)^{-1}\bigr]
\\
&&{}- \tr\bigl[\bigl\{\X'(\D+\la\I)^{-1}\X\bigr
\}^{-1}\X'(\D+\la\I)^{-1}\D\Q_{\W}\D(
\D+\la \I)^{-1}\X\bigr]
\\
&\geq& \frac{1}{ d_1+\la}\tr\bigl[\D^2\Q_{\W}\bigr] -
\frac{d_1 p }{ d_1+\la} \Ch_{\mathrm{max}}(\D\Q_{\W}),
\\
\Ch_{\mathrm{max}}\bigl(\D\Q_{\W}\D\A(\la)\bigr) &\leq&
\Ch_{\mathrm{max}}\bigl(\D\Q_{\W}\D(\D+\la\I)^{-1}\bigr) \leq
\frac{d_1}{d_1+\lah}\Ch_{\mathrm{max}}(\D\Q_{\W}),
\end{eqnarray*}
where $d_1\geq\cdots\geq d_k$.
Then,
\[
\tr[\D\Q_{\W}\D\A(\la)]/\Ch_{\mathrm{max}}\bigl(\D\Q_{\W}\D
\A(\la)\bigr)\geq\tr \bigl[\D^2\Q_{\W} \bigr]/\bigl\{
d_1 \Ch_{\mathrm{max}}(\D\Q_{\W})\bigr\} -p.
\]
Similarly, it is observed that
\begin{eqnarray*}
\tr\bigl[\D\Q_{\W}\D\A(\la)\bigr] &\geq& \frac{d_k}{ d_k+\la}\tr[\D
\Q_{\W}] - \frac{p }{ d_k+\la} \Ch_{\mathrm{max}}\bigl(\D^2
\Q_{\W}\bigr),
\\
\Ch_{\mathrm{max}}\bigl(\D\Q_{\W}\D\A(\la)\bigr) &\leq&
\frac{1}{ d_k+\lah} \Ch_{\mathrm{max}}\bigl(\D^2\Q_{\W}\bigr),
\end{eqnarray*}
which implied that $\tr[\D\Q_{\W}\D\A(\la)]/\Ch_{\mathrm{max}}(\D\Q_{\W}\D
\A(\la
))\geq d_k \tr[\D\Q_{\W}]/\break\Ch_{\mathrm{max}}(\D^2\Q_{\W}) -p$.
These provide the following sufficient condition.

\begin{prp}\label{prpFHmin2}
The constrained empirical Bayes estimator $\bmuh^{\mathrm{CEB}}(\lah,\t)$
improves on $\bmuh^{\mathrm{CM}}(\t)$ in terms of the conditional risk if the
following condition holds:
%
\begin{equation}
\label{eqnFHcond1} \max \biggl\{\frac{\tr[\D^2\Q_{\W}] }{ d_1 \Ch_{\mathrm{max}}(\D\Q_{\W})}, \frac{d_k \tr[\D
\Q_{\W}] }{\Ch_{\mathrm{max}}(\D^2\Q_{\W})} \biggr\}
\geq p+2 + \frac{k-p }{2}.
\end{equation}
\end{prp}

When $d_1=\cdots=d_k$ and $\Q=\I_m$, the condition (\ref{eqnFHcond1})
is written as $k-p \geq2(m+2)$, and improvement is guaranteed for
large $k$.
However, those sufficient conditions for the improvement are
restrictive in the case of different $d_i$'s with large $d_1$ and small $d_k$.

\textit{2. Unconditional risk.} We next investigate the dominance
property relative to the unconditional risk.
Let $\De^U=R^U(\pi,\bmuh^{\mathrm{CEB}}(\lah,\t)) - R^U(\pi,\bmuh^{\mathrm{CM}}(\t))$.
Since it is hard to evaluate $\De^U$ exactly, we shall approximate
$\De^U/k$ up to $\mathrm{O}(k^{-3/2})$ for large $k$.

\begin{prp}\label{prpFHmin3}
Assume that the elements of $\X$ and $\W$ are uniformly bounded and
$\X'\V^{-1}(\la)\*\X/k$ is positive definite and converges to a positive
definite matrix.
Assume also that $d_i$'s are bounded above and bounded away from zero.
Then, $\De^U/k$ is approximated as $\De^U/k=\De_{\mathrm{APR}}(\la)/k + \mathrm{O}(k^{-3/2})$,
where
%
\begin{eqnarray}
\label{eqnAPR} \De_{\mathrm{APR}}(\la) &=& - \tr\bigl[ \D\Q_{\W}\D
\V^{-1}(\la)\bigr]
+ \tr\bigl[ \bigl(\X'\V^{-1}(\la)\X
\bigr)^{-1}\X'\V^{-1}(\la)\D\Q_{\W}\D
\V^{-1}(\la)\X\bigr]\nonumber
\\[-4pt]\\[-12pt]
&&{}+ \tr\bigl[\D\Q_{\W}\D\V^{-3}(\la)\bigr]
\frac{2k }{(\tr[\V^{-1}(\la)])^2}.
\nonumber
\end{eqnarray}
A necessary condition for $\De_{\mathrm{APR}}(\la)\leq0$ is given by
%
\begin{equation}
\label{eqnFHcondn} \tr[ \D\Q_{\W}]\geq\tr\bigl[ \bigl(\X'
\D^{-1}\X\bigr)^{-1}\X'\Q_{\W}\X\bigr]
+ \tr\bigl[\Q_{\W}\D^{-1}\bigr] \frac{2 k}{(\tr[\D^{-1}])^2}.
\end{equation}
A sufficient condition for $\De_{\mathrm{APR}}(\la)\leq0$ is that
%
\begin{equation}
\label{eqnFHconds} \min_{\la>0} \biggl\{ \frac{\tr[\D\Q_{\W}\D\V^{-1}(\la)] }{\Ch_{\mathrm{max}}(\D\Q_{\W}
\D\V^{-1}(\la))} \biggr\}\geq p + 2
\frac{k \tr[\D^{-2}] }{(\tr[\D^{-1}])^2}.
\end{equation}
\end{prp}

The proof is given in the \hyperref[app]{Appendix}.
The approximation (\ref{eqnAPR}) was derived by Datta, Rao and Smith \cite{r8}.
When $\De_{\mathrm{APR}}(\la)\leq0$ for any $\la>0$, it is said that $\bmuh^{\mathrm{CEB}}(\lah,\t)$ improves on $\bmuh^{\mathrm{CM}}(\t)$ in terms of the
second-order approximation of the unconditional risk.
Using the same arguments as in (\ref{eqnFHcond1}), it follows from
(\ref{eqnFHconds}) that the inequality $\De_{\mathrm{APR}}(\la)\leq0$ holds if
%
\begin{equation}
\label{eqnFHcondss} \max \biggl\{ \frac{\tr[\D^2\Q_{\W}] }{ d_1 \Ch_{\mathrm{max}}(\D\Q_{\W})}, \frac{d_k \tr[\D
\Q_{\W}] }{\Ch_{\mathrm{max}}(\D^2\Q_{\W})} \biggr\}
\geq p + 2 \frac{k \tr[\D^{-2}]
}{(\tr[\D^{-1}])^2}.
\end{equation}
The necessary condition (\ref{eqnFHcondn}) is useful in the sense that
if the condition (\ref{eqnFHcondn}) is violated, then $\bmuh^{\mathrm{CEB}}(\lah
,\t)$ does not improve on $\bmuh^{\mathrm{CM}}(\t)$ in terms of the second-order
approximation of the unconditional risk.
This means that $\bmuh^{\mathrm{CEB}}(\lah,\t)$ should satisfy the condition
(\ref{eqnFHcondn}) at least.

\begin{remark}\label{remark1}
Propositions \ref{prpFHmin}, \ref{prpFHmin2} and \ref{prpFHmin3}
give us the conditions for the improvement by the constrained empirical
Bayes estimator (\ref{eqnCEB}).
By replacing $\Q_\W$ with $\Q$, these propositions can provide the
conditions under which the empirical Bayes estimator given in (\ref
{eqnEB}) improves on $\y$.
\end{remark}

\subsection{Unconstrained empirical Bayes estimator satisfying constraints}\label{secFH2}

In this subsection, we set up a prior distribution which results in an
unconstrained empirical Bayes and minimax estimator satisfying the
constraint in the Fay--Herriot model with heteroscedastic variances and
covariates as regressors.

\textit{Case 1}: $\t(\y)=\W'\y$.
Recall the arguments as in Case 1 of Section \ref{secnunc}.
Since $\bxi_3=(\H_1-\V_{12}\V_{22}^{-1}\H_2)\Q^{1/2}\bmu$ and we
set up
the linear regression structure $\X\bbe$ for $\bmu$, it may be
reasonable to assume the prior distribution $\bxi_3| \la\sim\Nc_{k-m}(\X_3\bbe, \la\I_{k-m})$ for $\bxi_3$ and to assume the uniform
prior for $\bxi_2$, where $\X_3=(\H_1-\V_{12}\V_{22}^{-1}\H_2)\Q^{1/2}\X
$, which is assumed to be of rank $p$.
Combining the contents in Sections \ref{secnunc} and \ref{secFH1}, we
get the empirical Bayes estimator given by
\[
\bmuh^{\mathrm{EB}1} = \Q^{-1/2}\H_1' \bigl\{
\bxih^{\mathrm{EB}}_3(\z_3) + \V_{12}
\V_{22}^{-1}\z_2\bigr\} + \Q^{-1/2}
\H_2'\z_2.
\]
Here the empirical Bayes estimator $\bxih^{\mathrm{EB}}_3(\z_3)$ is given as follows:
Note that $\z_3|\bxi_3 \sim\Nc_{k-m}(\bxi_3,\allowbreak \V_{11.2})$ and
$\bxi_3\sim
\Nc_{k-m}(\X_3\bbe, \la\I)$.
According to the arguments in Section \ref{secFH1}, we estimate $\la$
by $\lah=\max\{\la^*, 0\}$, where $\la^*$ is the solution of the
equation $\z_3'\A_3(\la^*)\z_3=k-m-p$ for $\A_3(\la) = \V_3^{-1}
- \X_3(\X_3'\V_3^{-1}\X_3)^{-1}\X_3'\V_3^{-1}$ for $\V_3=\V_{11.2}+\la\I$.
Then, the empirical Bayes estimator of $\bxi_3$ is written by
%
\begin{equation}
\label{eqnEB1} \bxih^{\mathrm{EB}}_3(\z_3) =
\z_3 - \V_{11.2}(\V_{11.2}+\lah\I_{k-m})^{-1}
\bigl\{\z_3-\X_3\bbeh_3(\lah) \bigr\}
\end{equation}
for $\bbeh_3(\la) = (\X_3'\V_3^{-1}\X_3)^{-1}\X_3'\V_3^{-1}\z_3$.

Clearly, $\bmuh^{\mathrm{EB}1}$ satisfies the constraint, namely, $\W'\bmuh^{\mathrm{EB}1}=\W'\y$.
Since $\Q_{\W}$ and $\D$ in Section~\ref{secFH1} correspond to $\I_{k-m}$ and $\V_{11.2}$, respectively.
The dominance results for $\bmuh^{\mathrm{EB}1}$ follow from Propositions \ref
{prpFHmin2} and \ref{prpFHmin3}.

\begin{prp}\label{prpEB1}
The unconstrained empirical Bayes estimator $\bmuh^{\mathrm{EB}1}$ satisfies the
constraint $\W'\bmuh^{\mathrm{EB}1}=\W'\y$.
It is also minimax in $\Ga$ in terms of the conditional risk if
%
\begin{equation}
\label{eqnFHas} \max \biggl\{ \frac{\tr[\V_{11.2}^2] }{\{\Ch_{\mathrm{max}}(\V_{11.2})\}^2 }, \frac{\Ch_{\mathrm{min}}(\V_{11.2}) }{\Ch_{\mathrm{max}}(\V_{11.2})}\tr[
\V_{11.2}] \biggr\} \geq p + 2 + \frac{k-p}{2}.
\end{equation}

In the sense of the second-order approximation relative to the
unconditional risk, a sufficient condition for $\bmuh^{\mathrm{EB}1}$ to improve
on $\y$ is
%
\begin{equation}
\label{eqnFHass} \max \biggl\{ \frac{\tr[\V_{11.2}^2] }{\{\Ch_{\mathrm{max}}(\V_{11.2})\}^2 }, \frac{\Ch_{\mathrm{min}}(\V_{11.2}) }{\Ch_{\mathrm{max}}(\V_{11.2})}\tr[
\V_{11.2}] \biggr\} \geq p + 2\frac{k\tr[\V_{11.2}^{-2}] }{(\tr[\V_{11.2}^{-1}])^2},
\end{equation}
and a necessary condition for the improvement is given by
%
\begin{equation}
\label{eqnFHan} \tr[ \V_{11.2}]\geq\tr\bigl[ \bigl(\X_3'
\V_{11.2}^{-1}\X_3\bigr)^{-1}
\X_3'\X_3\bigr] + 2 k/\tr\bigl[
\V_{11.2}^{-1}\bigr].
\end{equation}
\end{prp}

\textit{Case 2}: $\t(\y)=\t_0$.
Recall the arguments as in Case 2 of Section \ref{secnunc}.
If we assume the prior distribution that $\bxi_2=\bxi_0$ with
probability one for $\bxi_0=(\W'\Q^{-1/2}\H_2')^{-1}\t_0$, it is seen
that $\bxi_4=\bxi_1$.
Since $\bxi_1=\H_1 \Q^{1/2}\bmu$, it is reasonable to assume the prior
distribution $\bxi_1 | \la\sim\Nc_{k-m}(\H_1\Q^{1/2}\X\bbe, \la
\I_{k-m})$ for $\bxi_1$.
Then from (\ref{eqnGB2}), the generalized Bayes estimator is given by
%
\begin{equation}
\label{eqnEB2} \bmuh^{\mathrm{EB}2} = \Q^{-1/2}\H_1'
\bxih_1^{\mathrm{EB}}(\z_4) + \Q^{-1/2}
\H_2'\bxi_0,
\end{equation}
where $\bxih_1^{\mathrm{EB}}(\z_4)$ has the same form as $\bxih_3^{\mathrm{EB}}(\z_3)$
given in (\ref{eqnEB1}) except replacing $\z_3$ and $\X_3$ with $\z_4$
and $\H_1\Q^{1/2}\X$, respectively.
It is assumed that $\H_1\Q^{1/2}\X$ is of rank $p$.

Clearly, $\bmuh^{\mathrm{EB}2}$ satisfies the constraint, namely, $\W'\bmuh^{\mathrm{EB}2}=\t_0$.
The improvement of $\bmuh^{\mathrm{EB}2}$ follows from Propositions \ref
{prpFHmin} and \ref{prpEB1}.
When $\bmu$ is restricted to $\W'\bmu=\t_0$ or $\bxi_2=\bxi_0$,
$\bmuh^{\mathrm{EB}2}$ is minimax in $\Ga$ from Proposition \ref{prpprior2}.

\begin{prp}\label{prpEB2}
The unconstrained empirical Bayes estimator $\bmuh^{\mathrm{EB}2}$ satisfies the
constraint $\W'\bmuh^{\mathrm{EB}2}=\t_0$ and dominates the estimator $\Q^{-1/2}\H_1'\z_4 + \Q^{-1/2}\H_2'\bxi_0$ in terms of the conditional
risk if the condition $(\ref{eqnFHas})$ holds.
This implies the minimaxity of $\bmuh^{\mathrm{EB}2}$ within $\Ga_C$.
When $\bmu$ is restricted on $\W'\bmu=\t_0$, $\bmuh^{\mathrm{EB}2}$ is minimax
in $\Ga$ under the condition~$(\ref{eqnFHas})$.

In the sense of the second-order approximation relative to the
unconditional risk, a sufficient condition for $\bmuh^{\mathrm{EB}2}$ to improve
on $\y$ is given by $(\ref{eqnFHass})$, and a necessary condition for
the improvement is given by
%
\begin{equation}
\label{eqnFHa2n} \tr[ \V_{11.2}]\geq\tr\bigl[ \bigl(\X_4'
\V_{11.2}^{-1}\X_4\bigr)^{-1}
\X_4'\X_4\bigr] + 2 k/\tr\bigl[
\V_{11.2}^{-1}\bigr]
\end{equation}
for $\X_4=\H_1\Q^{1/2}\X$.
\end{prp}

\subsection{Simulation study}\label{secsim}

We investigate the unconditional risk behaviors of the constrained
estimators by simulation.
We consider the Fay--Herriot model (\ref{eqnEB}) with $k=15$, $\la
=1$ and four $d_i$-patterns: (a) 0.5, 0.5, 0.4, 0.3, 0.3; (b) 0.7, 0.6,
0.5, 0.4, 0.3; (c) 2.0, 0.6, 0.5, 0.4, 0.2; (d) 4.0, 0.6, 0.5, 0.4, 0.1.
Patterns (b)--(d) are treated by Datta, Rao and Smith \cite{r8}, and
pattern (a) is
less variable in $d_i$-values, while pattern (d) has larger variability.
There are five groups $G_1, \ldots, G_5$ and three small areas in each group.
The sampling variances $d_i$ are the same for areas within the same group.
For the matrix of covariates $\X$, the column vectors of $\X'$ are
generated as random vectors from $\Nc_p(\zero, (1-0.2)\I+0.2\j_k\j_k')$
where $\j_k$ is the $k$-dimensional vector with all the elements ones.
Each element of $\bbe$ is generated as $1+4u$ where $u\sim U(0,1)$, the
uniform distribution on $(0, 1)$.

\begin{table}
\tabcolsep=0pt
\caption{Values of unconditional risks of the constrained estimators
for $\la=1$}\label{tablerisk}
\begin{tabular*}{\textwidth}{@{\extracolsep{\fill}}lld{2.2}d{2.2}d{2.2}d{2.2}d{2.2}d{2.2}d{2.2}d{2.2}@{}}
\hline
&& & &\multicolumn{2}{l}{{Case 1}}&\multicolumn{2}{l}{{Case 2}}&
\multicolumn{2}{l@{}}{{Case 2$^*$}}\\[-5pt]
&& & &\multicolumn{2}{l}{\hrulefill}&\multicolumn{2}{l}{\hrulefill}&\multicolumn{2}{l@{}}{\hrulefill}\\
$\Q$ & $d_i$ & \multicolumn{1}{l}{$\y$} & \multicolumn{1}{l}{EB} & \multicolumn{1}{l}{CB} & \multicolumn{1}{l}{UC1} & \multicolumn{1}{l}{CB} & \multicolumn{1}{l}{UC2} & \multicolumn{1}{l}{CB}
& \multicolumn{1}{l}{UC2}\\
\hline
$\Q=\I$ & (a) & 6.00 & 4.76 & 4.84 & 4.88 & 9.24 & 9.40 & 4.48 & 4.64\\
& (b) & 7.51 & 5.53 & 5.63 & 5.70 & 9.68 & 9.90 & 5.22 & 5.45\\
& (c) &11.05 & 6.40 & 6.45 & 6.65 & 7.22 & 7.41 & 6.15 & 6.34\\
&{(d)}&16.88 & 6.60 & 6.61 & 7.26 &14.60 &19.92 & 6.47 &11.79\\
[5pt]
$\Q=\D^{-1}$&(a)&14.93 &11.92 &12.02 &12.28 &38.88 &39.14 &11.02 &11.28\\
&(b)&14.99 &11.49 &11.72 &11.87 &13.17 &13.31 &10.74 &10.89\\
&(c)&14.99 &10.90 &11.05 &11.76 &13.47 &14.17 &10.09 &10.80\\
&(d)&14.99 &10.57 &10.68 &12.07 &26.91 &28.30 & 9.70 &11.09\\
\hline
\end{tabular*}
\end{table}

In this simulation, we treat the case that $\W=\D^{-1}\j_k$,
$t_0=3\W'\X
\j_k$ and $\Q=\I, \D^{-1}$ for $m=1$ and $p=2$.
We compare the unconditional risks $R^U(\pi,\bmuh)$ for the five
estimators of $\bmu$: the crude estimator $\y$, the empirical Bayes
estimator EB given in (\ref{eqnEB}), the constrained empirical Bayes
estimator CB in (\ref{eqnCEB}), the unconstrained empirical Bayes
estimator UC1 in (\ref{eqnEB1}) for Case 1 and the unconstrained
empirical Bayes estimator UC2 in (\ref{eqnEB2}) for Case 2, where
Case 1 and Case~2 denote the constraints $\t(\y)=\W'\y$ and $\t
(\y
)=\t_0$, respectively.
The unconditional risks of these estimators are computed as average
values based on 10,000 simulation runs, and those values are reported
in Table \ref{tablerisk}, where Case 2$^*$ treats the unconditional
risks for $\bmu$ restricted to $\W'\bmu=t_0$.
It is noted that $\y$ and EB do not satisfy the constraints.
The values of the column of $\y$ correspond to the minimax risks for
Case 1 and Case 2$^*$, and it is revealed that EB, CB, UC1 and UC2
have smaller risks than $\y$.
For Case 1, the risks of the estimators CB and UC1 with the
constraints are slightly larger than those of $EB$.
It is interesting to note that the difference between Case 2 and
Case 2$^*$ supports Proposition \ref{prpmins}, namely, CB and UC2
improve on $\y$ when $\bmu$ is restricted to $\W'\bmu=t_0$, while their
maximum risks are beyond the risks of $\y$ without the restriction.

We next investigate whether the conditions for the improvement derived
in Sections \ref{secFH1} and~\ref{secFH2} are satisfied or not.
Table \ref{tablecond} reports this investigation where $+$ is marked if
the condition is satisfied, otherwise, $-$ is marked.
For improvement by CB, the sufficient condition relative to the
conditional risk is (\ref{eqnFHcond1}), denoted by SR, and the
sufficient and necessary conditions in terms of the second-order
approximation relative to the unconditional risk are given by (\ref
{eqnFHcondss}) and (\ref{eqnFHcondn}), respectively, denoted by
SR$^U$ and NR$^U$.
As noted in Remark \ref{remark1}, similar conditions for EB can be
given by SR, SR$^U$ and NR$^U$ by replacing $\Q_\W$ with $\Q$.
SR and SR$^U$ are given by (\ref{eqnFHas}) and (\ref{eqnFHass}) for
the improvement by UC1 and UC2.
The necessary conditions NR$^U$ for UC1 and UC2 are given by (\ref
{eqnFHan}) and (\ref{eqnFHa2n}), respectively.
As seen from Table \ref{tablecond}, the sufficient conditions SR under
the conditional risks for EB and CB are very restrictive in both cases
of $\Q=\I$ and $\Q=\D^{-1}$, and AR for UC1 and UC2 are also
restrictive for $\Q=\I$.
That is, the conditions SR are not satisfied in most cases.
It should be noted that this does not imply that those estimators do
not improve on $\y$, because the necessary conditions NR$^U$ are always
satisfied.
For the estimators UC1 and UC2, all the conditions for the improvement
are satisfied relative to the loss $(\bmuh-\bmu)'\D^{-1}(\bmuh-\bmu)$
for $\Q=\D^{-1}$.

\begin{table}
\tabcolsep=0pt
\caption{Whether the conditions for the improvement are satisfied or
not? When the condition is satisfied, $+$ is marked, and otherwise, $-$ is
marked}\label{tablecond}
\begin{tabular*}{\textwidth}{@{\extracolsep{\fill}}llllllllllllll@{}}
\hline
&& \multicolumn{3}{l}{{EB}}&\multicolumn{3}{l}{{CB}}
&\multicolumn{3}{l}{{UC1}}&\multicolumn{3}{l@{}}{{UC2}}\\[-5pt]
&& \multicolumn{3}{l}{\hrulefill}&\multicolumn{3}{l}{\hrulefill}&\multicolumn{3}{l}{\hrulefill}&\multicolumn{3}{l@{}}{\hrulefill}\\
$\Q$ & $d_i$ & SR & SR$^U$ & NR$^U$ & SR & SR$^U$ & NR$^U$ & SR & SR$^U$ &  NR$^U$ & SR & SR$^U$ & NR$^U$ \\
\hline
$\Q=\I$ & (a) & $-$ & $+$ & $+$ & $-$ & $+$ & $+$ & $-$ & $+$ & $+$ & $-$ & $+$ & $+$ \\
& (b)& $-$ & $+$ & $+$ & $-$ & $+$ & $+$ & $-$ & $+$ & $+$ & $-$ & $+$ & $+$ \\
&(c)& $-$ & $-$ & $+$ & $-$ & $-$ & $+$ & $-$ & $-$ & $+$ & $-$ & $-$ & $+$ \\
&(d)& $-$ & $-$ & $+$ & $-$ & $-$ & $+$ & $-$ & $-$ & $+$ & $-$ & $-$ & $+$ \\
[5pt]
$\Q= \D^{-1}$  &(a)& $+$ & $+$ & $+$ & $-$ & $+$ & $+$ & $+$ & $+$ & $+$ & $+$ & $+$ & $+$ \\
&(b)& $-$ & $+$ & $+$ & $-$ & $+$ & $+$ & $+$ & $+$ & $+$ & $+$ & $+$ & $+$\\
&(c)& $-$ & $-$ & $+$ & $-$ & $-$ & $+$ & $+$ & $+$ & $+$ & $+$ & $+$ & $+$ \\
&(d)& $-$ & $-$ & $+$ & $-$ & $-$ & $+$ & $+$ & $+$ & $+$ & $+$ & $+$ & $+$ \\
\hline
\end{tabular*}
\end{table}

\section{Concluding remarks}\label{secremark}

Benchmarking has been recognized as an important issue in small area
problems, and constrained Bayesian estimators have been studied in the
literature.
However, little has been known about decision-theoretic properties such
as admissibility and minimaxity for constrained generalized Bayes estimators.
In this paper, we have clarified admissibility, minimaxity and
dominance properties of constrained estimators by decomposing the
conditional risk function into two pieces: one depends on the
estimator, but the other does not depend on the estimator.
In the context of a multivariate normal population, we have provided a
canonical form, which allows us to establish admissibility and
inadmissibility of the constrained uniform-prior generalized Bayes estimator.
We have also derived a condition on the prior distribution such that
the resulting unconstrained generalized Bayes estimator automatically
satisfies the constraint.
Finally, we have provided constrained empirical Bayes and improved
estimators in the Fay--Herriot model.

Although a constrained empirical Bayes estimator is treated in Section
\ref{secFH}, it is not admissible. To develop admissible and minimax
estimators, we would need to consider hierarchical prior distributions
and to investigate admissibility and minimaxity of the resulting
hierarchical generalized Bayes estimators.
Berger and Robert \cite{r3}, Berger and Strawderman \cite{r4} and Kubokawa and Strawderman \cite{r15} have studied the admissibility and minimaxity of
hierarchical Bayes estimators.
The extension of their results to the setup of this paper seems a
reasonable goal and is one that we plan to study.

\begin{appendix}
\section*{Appendix: Proofs}\label{app}
\vspace*{-10pt}

\begin{pf*}{Proof of Proposition \protect\ref{prpFHmin}}
We first prove Proposition \ref{prpFHmin} which give us the sufficient
condition for the constrained empirical Bayes estimator $\bmuh^{\mathrm{CEB}}(\lah,\t)$ to improve on the constrained uniform-prior
generalized Bayes estimator $\bmuh^{\mathrm{CM}}(\t)$.
The arguments as in Shinozaki and Chang (\cite{r20,r21}) are useful for the proof.
The conditional risk difference of the two estimators is written as
\begin{eqnarray*}
\De&=& E\bigl[\bigl(\bmuh^{\mathrm{CEB}}(\lah,\t)-\bmu\bigr)'\Q
\bigl(\bmuh^{\mathrm{CEB}}(\lah,\t )-\bmu\bigr)|\bmu\bigr] - E\bigl[\bigl(
\bmuh^{\mathrm{CM}}(\t)-\bmu\bigr)'\Q\bigl(\bmuh^{\mathrm{CM}}(\t)-
\bmu\bigr)|\bmu\bigr]
\\
&=& E\bigl[\bigl(\bmuh^{\mathrm{EB}}(\lah)-\bmu\bigr)'
\Q_{\W}\bigl(\bmuh^{\mathrm{EB}}(\lah)-\bmu\bigr)|\bmu\bigr] - E\bigl[(
\y-\bmu)' \Q_{\W}(\y-\bmu)|\bmu\bigr],
\end{eqnarray*}
where $\Q_{\W}=\Q-\W(\W'\Q^{-1}\W)^{-1}\W'$.
It is noted that $\Q_{\W}$ is of rank $k-m$ and that $E[(\y-\bmu)'
\Q_{\W}(\y-\bmu)|\bmu]=\tr[\D\Q_{\W}]=\tr[\D\Q]-
\tr[\W'\D\W(\W'\Q^{-1}\W)^{-1}]$.
The risk difference is written as
%
\begin{equation}
\label{eqnA1} \De= -2E\bigl[(\y-\bmu)'\Q_{\W}\D\A(\lah)\y|
\bmu\bigr] + E\bigl[\y'\A(\lah )\D\Q_{\W}\D\A (\lah)\y|\bmu
\bigr].
\end{equation}
Using the Stein identity given in Stein \cite{r22}, we can rewrite the
cross product term as
\[
E\bigl[ (\y-\bmu)'\Q_{\W}\D\A(\lah) \y|\bmu\bigr] = E
\bigl[ \bnabla'\bigl\{\D\Q_{\W}\D\A(\lah)\y\bigr\} |\bmu
\bigr].
\]
Let $\G(\lah)=(g_{ij}(\lah)) = \D\Q_{\W}\D\A(\lah)$.
Then
\begin{eqnarray*}
\bnabla'\bigl\{ \G(\lah)\y\bigr\} &=& \sum
_{i,j} \frac{\partial}{\partial y_i} \bigl\{ g_{ij}(\lah)
y_j\bigr\}
\\
&=& \sum_i g_{ii}(\lah) + \sum
_{i,j} y_j \biggl\{ \frac{\dd}{\dd\la
}g_{ij}(
\la) \bigg|_{\la=\lah} \biggr\} \frac{\partial\lah}{
\partial y_i}
\\
&=& \tr\bigl[\D\Q_{\W}\D\A(\lah)\bigr] + \y' \biggl\{
\frac{\dd}{\dd\la}\A (\la) \bigg|_{\la
=\lah} \biggr\} \D\Q_{\W}\D(\bnabla
\lah),
\end{eqnarray*}
since $g_{ij}(\lah)$ depends on $\y$ through $\lah$.
Differentiating $\A(\la)$ with respect to $\la$ for $\A(\la)$
given in~(\ref{eqnEB}), we can see that
%
\begin{equation}
\label{eqnA2} \frac{\dd}{\dd\la} \A(\la) = - \A^2(\la),
\end{equation}
which can be used to get the expression $\De=E[\Deh|\bmu]$ where
\begin{eqnarray*}
\Deh(\lah) &=& -2 \tr\bigl[\D\Q_{\W}\D\A(\lah)\bigr] + 2
\y'\A^2(\lah)\D \Q_{\W}\D (\bnabla\lah)
\nonumber
\\
&&{}+ \y'\A(\lah)\D\Q_{\W}\D\A(\lah) \y
\end{eqnarray*}
for $\bnabla=(\partial/\partial y_1, \ldots, \partial/\partial y_k)'$.

Differentiating $\y'\A(\lah)\y=k-p$ with\vspace*{1pt} respect to $\y$ and using the
implicit function theorem, we get the equation $2\A(\lah)\y- \y'\A^2(\lah)\y\bnabla\lah=\zero$ in the case of $0< \lah$, or
\[
\bnabla\lah= \frac{2 }{\y'\A^2(\lah) \y} \A(\lah) \y I(0< \lah).
\]
Thus, $\Deh$ is expressed as
\begin{eqnarray*}
\Deh(\lah) &=& -2 \tr\bigl[\D\Q_{\W}\D\A(\lah)\bigr] + 4
\frac{\y'\A^2(\lah
)\D\Q_{\W}\D
\A(\lah)\y}{\y'\A^2(\lah)\y} I(0< \lah)
\nonumber
\\
&& {} + \y'\A(\lah)\D\Q_{\W}\D\A(\lah) \y,
\end{eqnarray*}
where $I(A)$ is the indicator function such that $I(A)=1$ if $A$ is
true, and otherwise, $I(A)=0$.
It is observed that
\begin{eqnarray*}
\y'\A(\lah)\D\Q_{\W}\D\A(\lah) \y &\leq& \y'
\A(\lah)\y\sup_{\mathbf{x}} \biggl\{ \frac{{\mathbf{x}}' \A(\lah)\D\Q_{\W}\D\A(\lah) \mathbf{x}}{{\mathbf{x}}'\A(\lah){\mathbf{x}}} \biggr\}
\\
&\leq&(k-p) \times\Ch_{\mathrm{max}}\bigl(\D\Q_{\W}\D\A(\lah)\bigr),
\\
\frac{\y'\A^2(\lah)\D\Q_{\W}\D\A(\lah)\y}{\y'\A^2(\lah)\y} I(\lah>0) &\leq& \sup_{\mathbf{x}} \biggl\{
\frac{{\mathbf{x}}'\A^2(\lah
)\D\Q_{\W}\D\A(\lah){\mathbf{x}}}{{\mathbf{x}}'\A^2(\lah){\mathbf{x}}} \biggr\}
\\
&=& \Ch_{\mathrm{max}}\bigl(\D\Q_{\W}\D\A(\lah)\bigr).
\end{eqnarray*}
Hence,
%
\begin{equation}
\label{eqnA3} \Deh(\lah) \leq-2 \tr\bigl[\D\Q_{\W}\D\A(\lah)\bigr] +
(k-p+4) \Ch_{\mathrm{max}}\bigl(\D\Q_{\W} \D\A(\lah)\bigr),
\end{equation}
which proves Proposition \ref{prpFHmin}.
\end{pf*}

\begin{pf*}{Proof of Proposition \ref{prpFHmin3}}
We next prove Proposition \ref{prpFHmin3}.
The unconditional risk difference can be written from (\ref{eqnA1}) as
\begin{eqnarray*}
\De^U &=& R^U\bigl(\pi,\bmuh^{\mathrm{CEB}}(\lah,\t)
\bigr) - R^U\bigl(\pi,\bmuh^{\mathrm{CM}}(\t)\bigr)
\\
&=& -2E\bigl[\bigl(\y-E[\bmu|\y]\bigr)'\Q_{\W}\D\A(\lah)
\y\bigr] + E\bigl[\y'\A(\lah)\D \Q_{\W}\D\A (\lah)\y\bigr].
\end{eqnarray*}
Noting that $E[\bmu|\y]=\y-\D\V(\la)^{-1}(\y-\X\bbe)$ and $\A
(\lah)\X
=\zero$, we see that
%
\begin{equation}
\De^U=-2E\bigl[ \u'\V(\la)^{-1}\D
\Q_{\W}\D\A(\lah)\u\bigr] + E\bigl[\u'\A (\lah)\D
\Q_{\W} \D\A(\lah)\u\bigr],
\end{equation}
where $\u$ is a random variable having $\Nc_k(\zero, \V(\la))$.
We shall derive the second order approximation of $\De^U/k$ up to $\mathrm{O}(k^{-1})$.
To this end, $\A(\lah)$ is approximated by the Taylor series
expansion as
\[
\A(\lah) = \A(\la) + \A^{(1)}(\la) (\lah-\la) + 2^{-1}
\A^{(2)}(\la) (\lah -\la)^2+\bigl[\mathrm{O}
\bigl(k^{-3/2}\bigr)\bigr]_{k\times k},
\]
where $\A^{(i)}(\la)=\partial^i \A(\la)/\partial\la^i$, $i=1,2$, and
$[\mathrm{O}(k^{-3/2})]_{k\times k}$ means that all elements of the matrix are
of $\mathrm{O}(k^{-3/2})$.
Then
\begin{eqnarray*}
&&E\bigl[\u'\A(\lah)\D\Q_{\W}\D\A(\lah)\u\bigr]
\\
&&\quad= E \bigl[ \u'\A(\la)\D\Q_{\W}\D\A(\la)\u+ 2
\u'\A(\la)\D\Q_{\W}\D\A^{(1)}(\la)\u(\lah-\la)
\\
&&\hspace*{10pt}\qquad{}+\u'\A(\la)\D\Q_{\W}\D
\A^{(2)}(\la)\u(\lah-\la)^2 + \u'
\A^{(1)}(\la)\D\Q_{\W}\D\A^{(1)}(\la)\u(\lah-
\la)^2 \bigr]\\
&&\hspace*{10pt}\qquad{}+ \mathrm{O}\bigl(k^{-1/2}\bigr)
\\
&&\quad= \tr\bigl[\D\Q_{\W}\D\A(\la)\bigr] - \tr\bigl[\bigl(
\X'\V^{-1}(\la)\X\bigr)^{-1}\X'
\V^{-1}(\la)\D\Q_{\W}\D\A(\la)\X\bigr]
\\
&&\qquad{}+ E \bigl[ 2\u'\V^{-1}(\la)\D\Q_{\W}
\D\A^{(1)}(\la)\u(\lah -\la) +\u'\V^{-1}(\la)\D
\Q_{\W}\D\A^{(2)}(\la)\u(\lah-\la)^2
\\
&&\hspace*{22pt}\qquad{}+ \u'\V^{-2}(\la)\D
\Q_{\W}\D\V^{-2}(\la)\u(\lah-\la)^2 \bigr] +
\mathrm{O}\bigl(k^{-1/2}\bigr),
\end{eqnarray*}
since $\u'\A(\la)\C\u=\u'\V^{-1}(\la)\C\u- \u'\V^{-1}(\la
)\X(\X'\V^{-1}(\la)\X)^{-1}\X'\V^{-1}(\la)\C\u=\break\u'\V^{-1}(\la)\C\u+
\mathrm{O}_p(1)$ for
a matrix $\C=[\mathrm{O}(1)]_{k\times k}$, and $\A^{(1)}(\la)=-\A^2(\la)$.
Similarly,
\begin{eqnarray*}
&&-2E\bigl[ \u'\V(\la)^{-1}\D\Q_{\W}\D\A(\lah)
\u\bigr]
\\
&&\quad= -2 \tr\bigl[ \D\Q_{\W}\D\A(\la)\bigr] -2E\bigl[
\u'\V^{-1}(\la)\D\Q_{\W}\D\A^{(1)}(\la)
\u(\lah-\la)\bigr]
\\
&&\qquad{}-E\bigl[ \u'\V^{-1}(\la)\D\Q_{\W}\D
\A^{(2)}(\la)\u(\lah-\la)^2\bigr] + \mathrm{O}
\bigl(k^{-1/2}\bigr).
\end{eqnarray*}
Since $E[\u'\V^{-2}(\la)\D\Q_{\W}\D\V^{-2}(\la)\u(\lah-\la
)^2]=E[\u'\V^{-2}(\la)\D\Q_{\W}\D\V^{-2}(\la)\u] \Var(\lah)+\mathrm{O}(k^{-1/2})$ and
$\A(\la
)\X=\zero$, it follows that
\[
\De^U = - \tr\bigl[ \D\Q_{\W}\D\A(\la)\bigr] + \tr\bigl[
\V^{-3}(\la)\D\Q_{\W}\D\bigr] \Var(\lah) + \mathrm{O}
\bigl(k^{-1/2}\bigr).
\]
It is noted that $\Var(\lah)=2k/(\tr[\V^{-1}(\la)])^2+\mathrm{O}(k^{-1/2})$ from
Datta, Rao and Smith \cite{r8}.
Hence, $\De^U/k$ can be approximated as $\De^U/k=\De_{\mathrm{APR}}(\la
)/k+\mathrm{O}(k^{-3/2})$, where\break $\De_{\mathrm{APR}}(\la)$ is given in (\ref{eqnAPR}).
A necessary condition for $\De_{\mathrm{APR}}(\la)\leq0$ is that $\De_{\mathrm{APR}}(0)\leq0$, which is given in (\ref{eqnFHcondn}).
To derive a sufficient condition, note that\break $\tr[ (\X'\V^{-1}(\la
)\X
)^{-1}\X'\V^{-1}(\la)\D\Q_{\W}\D\V^{-1}(\la)\X] \leq p \Ch_{\mathrm{max}}(\D\Q_{\W}
\D\V^{-1}(\la))$ and that $\tr[\D\Q_{\W}\D\V^{-3}(\la)]\leq\Ch_{\mathrm{max}}(\D
\Q_{\W}\D\V^{-1}(\la))\tr[\V^{-2}(\la)]$.
By making the differentiation, it can be verified that $\tr[\V^{-2}(\la
)]/(\tr[\V^{-1}(\la)])^2$ is decreasing in $\la$, so that $\tr[\V^{-2}(\la)]/(\tr[\V^{-1}(\la)])^2\leq\tr[\D^{-2}]/(\tr[\D^{-1}])^2$.
Thus,
\[
\De_{\mathrm{APR}}(\la) \leq - \tr\bigl[ \D\Q_{\W}\D\V^{-1}(
\la)\bigr] + \Ch_{\mathrm{max}}\bigl(\D\Q_{\W}\D\V^{-1}(\la)
\bigr) \bigl\{ p + 2 k \tr\bigl[\D^{-2}\bigr]/\bigl(\tr\bigl[
\D^{-1}\bigr]\bigr)^2 \bigr\},
\]
which is expressed as (\ref{eqnFHconds}).
Therefore, we get Proposition \ref{prpFHmin3}.
\end{pf*}
\end{appendix}

\section*{Acknowledgements}

We are grateful to the Editor, the Associate Editor and the reviewers
for their valuable comments and helpful suggestions.
This research was supported in part by Grant-in-Aid for Scientific
Research $\#$ 21540114 and $\#$ 23243039 from Japan Society for the
Promotion of Science.
This work was partially supported by a grant from the Simons Foundation
($\#$ 209035 to William Strawderman).



%

\printhistory

\end{document}